\documentclass[secnum]{rims29}
\usepackage{amssymb, amsmath}
\usepackage[dvipdfm]{hyperref}
\usepackage{enumerate}
\usepackage{epic,eepic}
\usepackage{graphicx}
\usepackage{float}
\usepackage[mathscr]{eucal}
\newtheorem{crl}{Corollary}[section]

\newtheorem{lmm}{Lemma}[section]

\newtheorem{prp}{Proposition}[section]

\newtheorem{thm}[lmm]{Theorem}

\newtheorem{prop}[lmm]{Proposition}
\newtheorem{fact}[lmm]{Fact}
\newtheorem{prob}[lmm]{Problem}
\theoremstyle{definition}
\newtheorem{dfn}{Definition}[section]
\newtheorem{exa}{Example}
\theoremstyle{remark}
\newtheorem{remark}[lmm]{Remark}
\newcommand{\Fka}{\mathcal{F}_{k,a}}%
\def\Sol{{\mathcal{S}}ol}
\def \set#1#2{\{{#1}:{#2}\}}
\def \indefriem#1#2{X_{#1}^{#2}}
\title{Algebraic analysis of minimal representations}
\dedicatory{Dedicated to Mikio Sato
 whose pioneering work in 
\\ algebraic analysis
 has been an inspiration for me.  }           
\author{Toshiyuki \textsc{Kobayashi}
          \footnote{Graduate School of Mathematical Sciences, The University of
 Tokyo,
3-8-1 Komaba, Meguro, Tokyo, 153-8914 Japan.}}
\AuthorHead{Toshiyuki Kobayashi}         
\support{Partially supported by Grant-in-Aid
      for Scientific Research (B) (18340037), Japan
      Society for the Promotion of Science}
\VolumeNo{4x}           
\YearNo{200x}           
\PagesNo{000--000}      
\communication{JJ, December 31, 2009;
        Revised *.}      
\begin{document}
%

\maketitle


\begin{abstract}      
Small representations of a group bring us
 to large symmetries 
 in a representation space.  
Analysis on minimal representations
 utilises large symmetries
 in their geometric models, 
 and serves as a driving force in creating new interesting problems
 that interact with other branches 
 of mathematics.  

This article discusses the following three topics
 that arise from minimal representations
 of the indefinite orthogonal group:
\begin{enumerate}
\item[1.]
construction of conservative quantities for ultra-hyperbolic equations,
\item[2.]
quantative discrete branching laws,
\item[3.]
deformation of the Fourier transform
\end{enumerate}
 with emphasis on the prominent roles
 of Sato's idea on algebraic analysis.  
\end{abstract}

\section{Introduction}
\label{sec:1}

The aim of this article is to highlight the prominent roles of Sato's
idea on hyperfunctions and $\mathcal{D}$-modules
in the new developments on 
analysis of minimal representations \cite{xbko,xkmano2,xkmanoAMS,xkors3}.

Minimal representations are the simplest, infinite dimensional
`unipotent representations'.
They are building blocks of unitary
representations.
Segal--Shale--Weil representation is a classic example for the split
simple group of type $C$.
There has been an active study on minimal representations of reductive groups,
mostly by \textit{algebraic approaches} since 1990s both in the real and
in the $p$-adic fields
\cite{xBiZi,xbk,xgansavin,xGrWal,xKa,xkazsav,xkohcrcras,xKo,xMo,xSab,xTo,xV,xHuZu}.

On the other hand,
I believe that 
 \textit{geometric analysis} on minimal
representations is also a promising area,
and have been advocating its study 
based on the following change of viewpoints:
\begin{align}\label{eqn:philo}
&\text{\textbf{small} representations of a group}
\nonumber
\\
={}
&\text{\textbf{large} symmetries in a representation space}.
\end{align}
The terminology `minimal representations' is defined 
inside representation theory
(i.e.\ the annihilator in the universal enveloping algebra is the
Joseph ideal, see e.g.\ \cite{xgansavin}),
and the corresponding `largest symmetries' 
are expected to serve as a driving force in creating new
interesting areas of mathematics
 even outside representation theory.

The `largest symmetries' in representation spaces of minimal representations
 may be also observed in branching laws.
Indeed, as we shall see in Theorem \ref{thm:twistpull}, 
it may well happen
 that \textit{broken symmetries} of minimal representations
 reduce to analysis on certain
semisimple symmetric spaces (see also \cite{xkcheck,xkors2, KOP}).
This observation indicates that
analysis on minimal representation involves higher
symmetries than those for (traditional) analysis on symmetric spaces.

We focus on the minimal representation of a simple group
 of type D.  
This is just a single irreducible representation, 
 however,
 it turns out that geometric analysis on its various models is surprisingly rich.  
Indeed, 
papers devoted
 to this single representation
 in very recent years already
 exceed 500 pages, 
 giving rise to the interactions with the following topics:

\begin{itemize}
\item[$\bullet$]
conformal geometry 
 for general pseudo-Riemannian manifolds
 \cite{xkcheck, xkors1}, 
\item[$\bullet$]
Dolbeault cohomologies on open complex manifolds \cite{xkohcrcras}.
\item[$\bullet$]
conservative quantities for ultra-hyperbolic equations 
 \cite{xkors3}, 
\item[$\bullet$]
breaking symmetries and discrete branching laws
 \cite{xkors2}, 
\item[$\bullet$]
 Schr\"odinger model and the unitary inversion operator
 \cite{xkmanoprojp, xkmanoAMS}, 
\item[$\bullet$]
deformation of Fourier transforms
 \cite{xbko}, 
\item[$\bullet$]
holomorphic semigroup 
 \cite{xkmano1, xkmano2}, 
\item[$\bullet$]
new special function theory
 for fourth order differential operators
 \cite{xhkmm1, xhkmm2}.  
\end{itemize}

In this article, 
 we choose three topics among them, 
 and try to explain their flavours
 in Sections \ref{sec:2},
 \ref{sec:3}
 and \ref{sec:4}, 
 respectively with emphasis on the role of 
Sato's idea on algebraic analysis, 
 both in philosophy
 and in techniques.      
For the reader's convenience,
 we list some representation theoretic properties
 of our minimal representation
 in the Appendix.  

\section{Conservative quantities
 for ${\mathcal {D}}$-modules. }
\label{sec:2}

The energy of a wave is a conservative quantity 
 for the wave equation, 
 namely, 
 it is invariant under time-translations.  
In this section,
 we discuss higher symmetries coming from conformal transformations.  
By using the idea of Sato's hyperfunctions \cite{xkkk, xsato}, 
 we construct conservative quantities
 for specific ultra-hyperbolic equations
 (see Theorem \ref{thm:definpro}).

\subsection{Yamabe operator and conformal geometry}
\label{subsec:2.1}

A diffeomorphism $h$ of a Riemannian manifold $(X,g)$ is said to be {\it{conformal}}
 if there exists a positive-valued function 
 $\Omega(h, \cdot)$ on $X$
 such that 
$$
   h^{\ast} g_{hx} = \Omega(h,x)^2 g_x
  \quad
  \text{for}\quad
  x \in X.  
$$
It is isometry
 if $\Omega(h,\cdot) \equiv 1$.  
We write 
$$
  {\operatorname{Isom}}(X,g) \subset \operatorname{Conf} (X,g)
$$
for the groups consisting of isometries 
 and conformal diffeomorphisms, 
 respectively.  
The same notation will be applied to a more general setting 
 where $g$ is a  non-degenerate symmetric tensor,
 namely, 
 to an indefinite-Riemannian manifold.

The invariance for the Laplacian $\Delta_X$ characterizes isometries among diffeomorphisms of $X$.  
In other words, 
 a non-isometric transformation on $(X,g)$ does not preserve $\Delta_X$.  
However,
 the Laplacian $\Delta_X$ is still subject to the following covariance
 under conformal transformations:
\begin{equation}
\label{eqn:crep}
  \varpi_{\frac{n+2}{2}}(h) \circ \widetilde{\Delta}_X
  =
  \widetilde{\Delta}_X \circ \varpi_{\frac{n-2}{2}}(h)
  \quad
  \text{for any }\,\,
  h \in \operatorname{Conf}(X,g), 
\end{equation}
where $n$ is the dimension of $X$, 
 $\operatorname{Scal}_X$ is the scalar curvature, 
 and 
\begin{alignat*}{2}
 \widetilde{\Delta}_X :=& \Delta_X - \frac{n-2}{4(n-1)} \operatorname{Scal}_X
 \quad
  &&\text{(the Yamabe operator), }
\\
 \varpi_{\lambda} (h) f (x) := & \Omega(h^{-1}, x)^{\lambda} f(h^{-1}  x)
  \quad
  &&\text{for $f \in C^{\infty}(X)$.} 
\end{alignat*}

The formula \eqref{eqn:crep} implies 
 that the operator $\Delta_X$
 (or $\widetilde{\Delta}_X$)
 is not conformally invariant,  
 but the ${\mathcal {D}}$-module generated by 
 $\widetilde{\Delta}_X$ is conformally invariant !  
As far as the solutions
 are concerned, 
 the invariance of the ${\mathcal{D}}$-module is sufficient.  
Namely, 
 by putting 
\begin{equation}
\label{eqn:YSol}
  \Sol(\widetilde{\Delta}_X)
  :=
  \set{f \in C^{\infty}(X)}
       {{\Delta_X} f =  \operatorname{Scal}_X f}, 
\end{equation}
 we get: 
 
\begin{fact}
\label{fact:confrep}
The conformal group $\operatorname{Conf}(X,g)$
 preserves $\Sol(\widetilde{\Delta}_X)$ via $\varpi_{\frac{n-2}{2}}$.
\end{fact} 

See \cite[Theorem 2.5]{xkors1} for the proof.

\begin{remark}
{\rm{
\begin{enumerate}
\item[{\rm{1)}}]
Other eigenspaces
 $\Sol(\widetilde{\Delta}_X - \lambda)$ are not conformally invariant if $\lambda \ne 0$.  
\item[{\rm{2)}}]
 It may be better formulated  
 if we use the ring of twisted differential operators
 acting on sections of the line bundle ${\mathcal{L}}_{\frac{n-2}{2}}$.  
\item[{\rm{3)}}]
The differential equation 
 $\widetilde{\Delta}_X f=0$, 
 namely, 
 $\Delta_X f={\operatorname{Scal}}_X f$, 
 is elliptic, hyperbolic, or ultra-hyperbolic, respectively,
 if $(X,g)$ is Riemannian, Lorentzian, or of general signature, 
 respectively.  
\end{enumerate}
}}
\end{remark}

Then a general problem is:
\begin{prob}[
{\rm{(see {\cite[Problem C]{xkcheck}})}}]
\label{prob:innersol}
\begin{enumerate}
\item[{\rm{1)}}]
Does there exist an invariant inner product on an appropriate subspace of $\Sol(\widetilde{\Delta}_X)$?  
\item[{\rm{2)}}]
If yes,
 construct it explicitly.  
\end{enumerate}
\end{prob}

Such an inner product may be seen as a conservative quantity
 for the solution to the equation
$
  \widetilde{\Delta}_X f =0.  
$
Problem \ref{prob:innersol} does not find a final answer 
 in the general setting.  
We shall give a partial answer in the flat case
 (see Theorem \ref{thm:definpro} below).

\subsection{Conservative quantities}
\label{subsec:2.2}

Let ${\mathbb{R}}^{p,q}$ be the Euclidean space ${\mathbb{R}}^{p+q}$
 endowed with the flat indefinite-Riemannian structure
$$
   ds^2=dx_1^2 + \cdots + dx_p^2 - dx_{p+1}^2
   - \cdots - dx_{p+q}^2.
$$
Then, 
  the corresponding Laplace--Beltrami operator takes the form:
$$
  {\square}_{p,q}
  :=
  \frac {{\partial}^2}{\partial x_{1}^2}
  +
  \cdots
  +
  \frac {{\partial}^2}{\partial x_{p}^2}
  -
  \frac {{\partial}^2}{\partial x_{p+1}^2}
  -
  \cdots
  -
  \frac {{\partial}^2}{\partial x_{p+q}^2}.  
$$
Obviously, 
 the scalar curvature on ${\mathbb{R}}^{p,q}$
 vanishes identically.  
Hence, the Yamabe operator on ${\mathbb{R}}^{p,q}$
 coincides with $\square_{p,q}$.  
The space of solutions to $\square_{p,q} f=0$,
 denoted by $\Sol(\square_{p,q})$,
 is obviously invariant
 under the motion group
$$
     \operatorname{Isom}({\mathbb{R}}^{p,q})
    \simeq O(p,q) \ltimes {\mathbb{R}}^{p+q}.   
$$
It was proved in \cite[Theorem 4.7]{xkors3}
 that $\Sol(\square_{p,q})$ has even larger symmetries
  if $p+q$ is even, 
 namely,
 by the indefinite orthogonal group
$$
   G:= O(p+1, q+1)
      =\{g \in GL(p+q+2, {\mathbb{R}}):
         {}^{t\!} g I_{p+1,q+1} g = I_{p+1, q+1}\}
$$
 acting on ${\mathbb{R}}^{p+q}$
 as M\"obius transforms.  
(To be more precise, 
 $G$ preserves the space
 $\Sol_0(\square_{p,q})$ of smooth solutions
 with certain decay conditions at infinity
 together with their derivatives.)

\begin{remark}
\begin{enumerate}
\item[1)]
The parity condition on $p+q$ is crucial.
In fact, a theorem of Howe and Vogan \cite{xV} asserts that there does not exist an infinite dimensional representation of $G$ whose Gelfand--Kirillov dimension is $p+q-1$ if $p+q$ is odd and $p,q>3$.
\item[2)]
$\Sol_0 (\square_{p,q})$ is defined as the twisted pull-back of smooth functions on the conformal compactification of $\mathbb{R}^{p+q}$.
See \cite{xkors3} for details.
\end{enumerate}
\end{remark}

Problem \ref{prob:innersol}
 in this specific setting is stated as:  

\begin{prob}
\label{prob:cq}
Find a $G$-invariant inner product 
 on $\Sol_0(\square_{p,q})$
 if exists.  
\end{prob}

\subsection{Unitarizability versus unitarization.}
\label{subsec:2.3}

If $p,q >0$ and $p+q$ is even and greater than two, 
 then we can tell a priori 
 that the representation on $\Sol_0(\square_{p,q})$
 is unitarizable and irreducible
 (e.g. \cite{xBiZi, xkors1})
  by algebraic techniques.  
Namely, 
 we know the existence and the uniqueness of a $G$-invariant inner product
 on $\Sol_0(\square_{p,q})$
 in this case.

What we seek for in Problem \ref{prob:cq}
 is not merely an abstract {\it{unitarizability}}
 but the {\it{unitarization}} of the representation space
 for a concrete geometric model, 
 namely, 
 the construction of the invariant inner product.  
Then,
 there are two approaches to the unitarization --- one is easier, 
 and the other is more challenging
 as discussed below.

An easier approach to Problem \ref{prob:cq} is 
 to write the inner product
 by using the integral representation of solutions.  
Such an integral formula was given in \cite[Theorem 4.7]{xkors3} by using an explicit formula of the Green kernel \cite{xkcheck, xkors3}.   
The disadvantage of this approach is that the formula
 of the inner product involves
 a preimage of the integral representation, 
 which is not canonically given.

A second approach is to use
 an expansion of solutions
 into a countable sum of better understood solutions, 
 and then to give a Parseval--Plancherel type theorem.  
We shall discuss this approach in Section \ref{sec:3}.

A more intrinsic approach is to find a formula of the inner product directly
 without using the integral representation of solutions. 
A hint to this is the well-known formula of energy
 to the wave equation,   
which is given by the integration
 of the Cauchy data on the hyperplane 
 ($t=\text{constant}$)
 in the space-time, 
 see \eqref{eqn:E}.  
(We note that, 
 however, 
 the energy is not conformally invariant
 but invariant only under time-translations.)

In order to explain the second approach, let us set up some notation.
We recall that any non-characteristic hyperplane in ${\mathbb{R}}^{p,q}$
 is written as 
\begin{equation}
\label{eqn:avc}
  \alpha \equiv \alpha_{v,c} :=
  \{x \in {\mathbb{R}}^{p+q}:
    (x,v)_{{\mathbb{R}}^{p,q}}=c\}
\end{equation}
 for some $c \in {\mathbb{R}}$ and $v \in {\mathbb{R}}^{p,q}$ such that 
 $(v,v)_{{\mathbb{R}}^{p,q}}= \pm 1$.  
Fix such $v$, 
 and express a function $f$ on ${\mathbb{R}}^{p+q}$ as Sato's hyperfunction 
 (\cite{xsato})
 in the direction of $v$,
 namely,  
\begin{equation}
\label{eqn:hyperf}
     f(x)= \lim_{\varepsilon \downarrow 0} (f_+ (x + \sqrt{-1} \varepsilon v) - f_-(x - \sqrt{-1} \varepsilon v)).  
\end{equation}
Here, 
 $f_{\pm}(x+tv)$ is a holomorphic function
 of one variable $t$ near the real axis in $\pm \operatorname{Im}t>0$.  

We set 
$$
     \frac{\partial f_{\pm}}{\partial \nu}(x):= 
     \frac{\partial}{\partial t}|_{t=0} f(x + t v)
    \qquad\text{(normal derivative), }
$$
 and introduce a function $Q_{\alpha} f$ on the hyperplane $\alpha \equiv \alpha_{v,c}$
 by 
$$
    Q_{\alpha} f :=\frac{1}{\sqrt{-1}}
         (f_+ \overline{\dfrac{\partial f_+}{\partial\nu}}
    -f_- \overline{\dfrac{\partial f_-}{\partial \nu}}).  
$$  
Finally,
 we define
\begin{equation}
\label{eqn:Qff}
   (f,f):= \int_{\alpha} Q_{\alpha} f.   
\end{equation}

The right-hand side of \eqref{eqn:Qff} does not always converge, 
 but it makes sense 
 if $f$ satisfies suitable decay conditions, 
 say, 
 $f \in \Sol_0(\square_{p,q})$.  
Then, 
 we can give an answer to Problem \ref{prob:cq} as follows:
\begin{thm}[{\rm{(see {\cite[Theorem 6.2]{xkors3}}, also \cite{xkcheck})}}]
\label{thm:definpro}
\begin{enumerate}
\item[{\rm{1)}}]
For $f \in \Sol_0(\square_{p,q})$, 
\eqref{eqn:Qff} is independent of the choice 
 of the pair $(f_+, f_-)$ in the expression \eqref{eqn:hyperf} and of the hyperplane $\alpha$.  
\item[{\rm{2)}}]
$(f,f) \ge 0$
 for any $f \in \Sol_0(\square_{p,q})$.  
The equality holds if and only if $f=0$.  
\item[{\rm{3)}}]
The polarization of the norm \eqref{eqn:Qff} yields a $G$-invariant inner product on $\Sol_0(\square_{p,q})$.  
\end{enumerate}
\end{thm}

We denote by $\overline{\Sol_0(\square_{p,q})}$
 for the Hilbert space obtained as the completion of $\Sol_0(\square_{p,q})$.  
Then, 
 we get a unitary representation
 of $G=O(p+1,q+1)$, 
 to be denoted by $\varpi\equiv \varpi^{p+1, q+1}$, 
 on $\overline{\Sol_0(\square_{p,q})}$.  
 It turns out that this is irreducible and
 a minimal representation of $G$.  
See Section \ref{sec:5}
 for representation theoretic properties of $\varpi$.  

\begin{remark}
\label{rem:definpro}
{\rm{
The assertion 1) in Theorem \ref{thm:definpro}
 is a part of the invariance
 of the inner product $(\, , \, )$
 because any non-characteristic hyperplane is conjugate
 to either $x_1=0$ or $x_{p+q}=0$
 by the motion group
$
   \operatorname{Isom}({\mathbb{R}}^{p,q}) \simeq O(p,q) \ltimes {\mathbb{R}}^{p+q}.  
$
We note that $G$
 contains $\operatorname{Isom}({\mathbb{R}}^{p,q})$
 as a proper subgroup.  
}}
\end{remark}

The proof of Theorem \ref{thm:definpro}
 was given in \cite{xkors3}
 by using some representation theoretic results
 of the representation $\varpi$.  
It might be interesting to find a proof that does not depend on group theory but only on geometry such as Stokes' theorem.
We then pin down this as 
 an open problem:  

\begin{prob}
\label{prob:gp}
Give a purely geometric proof to Theorem \ref{thm:definpro}.  
\end{prob}

\subsection{Energy generator}
\label{subsec:2.4}  

Our conformally invariant inner product \eqref{eqn:Qff}
 is very close to the energy of the wave,
 where one integrates Cauchy data 
 on the zero time hyperplane. 
We end this section
 by making more explicit its connection.

For $p=1$, 
 let us introduce time and space coordinates
 $(t;x)$
 instead of the previous coordinates
 $(x_1, \cdots, x_p;x_{p+1}, \cdots, x_{p+q})$. 
Then, 
 the energy of the wave $f$
 is given by
\begin{equation}
\label{eqn:E}
  {\mathcal {E}}(f)
  =
  \frac{1}{2}\int_{{\mathbb{R}}^q}
  (|f_t|^2+|\nabla f|^2)dx.  
\end{equation}
Then, 
 in terms of the inner product 
 \eqref{eqn:Qff}, 
 ${\mathcal {E}}(f)$ is written as 
$$
    (f,|H|f)=(f^+,Hf^+)-(f^-,H f^-)
$$
 where $H=i \partial_t$
 is the energy generator
 (infinitesimal time translations).  
Since the energy generator $H$ is invariant 
 under time-translations (i.e. invariant by a one-dimensional subgroup of $G$)
 and the inner product $(\, , \, )$ is invariant under the whole group $G$, 
 ${\mathcal{E}}(f)$ is also invariant under time-translations.
This explains the classical fact that the energy
 \eqref{eqn:E}
 is a conservative quantity in the narrow sense that it is 
 independent of which constant-time
 hyperplane
 we integrate over.

\section{Quantative branching laws}
\label{sec:3}

In Section \ref{sec:1}, 
 we have given a concrete formula of the conformally invariant inner product
 on the minimal representation
 ({\it{conservative quantities}}).  
It is given by the integral over hyperplanes.  
Yet another formula of the same inner product
 will be given as a countable sum
 of well-understood quantities.

This is a Parseval-type theorem
 (see Theorem \ref{thm:twistpull}), 
 which is built on a \lq{}good expansion theorem' of solutions.  
Such an expansion theorem is obtained
 as a special case of the general theory
 of discretely decomposable restrictions
 of unitary representations
 (see Theorem \ref{thm:deco}).  
We will see that algebraic analysis provides a powerful method
 to branching problems in representation theory
 (cf. Problem \ref{prob:deco} below).

\subsection{Breaking symmetries and discrete decomposability}
\label{subsec:3.1}
Suppose $\pi: G \to GL({\mathcal{H}})$
 is a unitary representation
 of a Lie group $G$.  
Given a subgroup $G'$ of $G$, 
 and consider the broken symmetry, 
 namely, 
 the restriction $\pi|_{G'}$.   
In general, 
 the restriction $\pi|_{G'}$ decomposes into a direct integral of 
 irreducible representations
 of $G'$.  
Our concern here is with:
\begin{prob}
[{\rm{(see \cite{xk:1, xk:2})}}]
\label{prob:deco}
For which triple $(G, G',\pi)$
 does the restriction $\pi|_{G'}$ decompose discretely
 with finite multiplicities?
\end{prob}

It often happens that the irreducible decomposition
 of the restriction $\pi|_{G'}$
 ({\it{branching law}}) contains continuous spectrum
 if $G'$ is non-compact.  
Even worse,
 each irreducible representation
 of $G'$ may occur in the branching law with infinite multiplicities.  
Thus,
 Problem \ref{prob:deco} seeks for a very nice class of branching laws.

Now, 
 let us fix some notation for a real reductive group $G$.  
Let $K$ be a maximal compact subgroup of $G$,
 $T$ a maximal torus of $K$, 
 and ${\mathfrak {t}}$, ${\mathfrak {k}}$ the Lie algebras
 of $T$, $K$, 
 respectively.  
We choose the set $\Delta^+({\mathfrak {k}}, {\mathfrak {t}})$
 of positive roots,
 and denote 
 the dominant Weyl chamber by ${\mathfrak {t}}_+ (\subset \sqrt{-1}{\mathfrak {t}}^*)$.  
We also fix a $K$-invariant inner product on ${\mathfrak {k}}$, 
 and regard $\sqrt{-1} {\mathfrak {t}}^*$ as a subset
 of $\sqrt{-1}{\mathfrak {k}}^*$.

Suppose that $K'$ is a closed subgroup of $K$.  
The group $K$ acts on the homogeneous space $K/K'$
 from the left, 
 and then on the cotangent bundle $T^{\ast}(K/K')$
 in a Hamiltonian fashion.   
 We write
$$
  \mu: T^{\ast}(K/K') \to \sqrt{-1}{\mathfrak {k}}^*
$$
for the momentum map, 
 and define the following closed cone by 
$$
    C_K(K') := \operatorname{Image}\mu \cap {\mathfrak {t}}_+.  
$$
For a subgroup $G'$ of $G$, 
 we shall consider $C_K(K')$
 by setting $K' := K \cap G'$.

Next, 
 suppose that $\pi$ is a (reducible) representation 
 of a compact Lie group $K$.  
The asymptotic $K$-support of $\pi$,
 to be denoted by ${\operatorname{AS}}_K(\pi)$,
 was introduced by Kashiwara and Vergne \cite{xKaVe}
 as the asymptotic cone of the $K$-types of $\pi$.  
From definition ${\operatorname{AS}}_K(\pi)=\{0\}$
 if $\dim \pi < \infty$.  
For a representation $\pi$ of $G$, 
 we can define $
{\operatorname{AS}}_K(\pi)$ 
 for the restriction $\pi|_K$.

We are ready to state an answer 
 to Problem \ref{prob:deco}:
\begin{thm}
[{\rm{(see \cite{xk:3})}}]
\label{thm:deco}
Suppose that $\pi$ is a unitary representation of $G$
 of finite length, 
 and that $G'$ is a closed subgroup of $G$.  
We set $K'= K \cap G'$.  
If 
\begin{equation}
\label{eqn:criterion}
C_K(K') \cap {\operatorname{
AS}}_K(\pi) = \{0\},
\end{equation}
then the restriction $\pi|_{G'}$ decomposes discretely into a direct sum 
of irreducible unitary representations of $G'$
 with finite multiplicities.  
\end{thm}
 
An upper estimate of the singularity spectrum
 of the hyperfunction character of $\pi$ plays a crucial role
 in the proof of Theorem \ref{thm:deco}. 
In particular, 
 the assumption \eqref{eqn:criterion} assures
\begin{equation}
\label{eqn:TR}
\text{Restriction and Trace ({\it{hyperfunction character}}) commute.  }
\end{equation}
Here, we remark that the character of an infinite dimensional representation $\pi$,
\[
\operatorname{Trace} \pi(g) \qquad (g \in G)
\]
does not make sense as an ordinary function
 because $\operatorname{Trace} \pi(e) = \dim \pi = \infty$.
Harish-Chandra proved that $\operatorname{Trace} \pi$ is well-defined
 as a distribution on $G$
 if $\pi$ is an irreducible unitary representation of a real reductive group $G$,
 and proved further that $\operatorname{Trace} \pi$ belongs to $L^1_{\mathrm{loc}}(G)$.
On the other hand, the restriction $\operatorname{Trace} \pi|_K$ is not locally integrable on $K$ any more
(see Atiyah \cite{xAtiyah}).
What \eqref{eqn:TR} means is that 
\[
\operatorname{Trace} (\pi|_{K'}) = \operatorname{Trace} (\pi)|_{K'}
\]
as an identity of hyperfunctions (or distributions) on $K'$.
See \cite[Theorem 2.8]{xk:3} for the proof.  
We also refer to the lecture notes \cite{xkprog05} for heuristic ideas of the proof.

Recently, 
 Hansen, Hilgert, and Keliny \cite{xhhdistr}
 has given an alternative proof of Theorem \ref{thm:deco}
 by replacing Sato's hyperfunctions with Schwartz's distributions.  
See also \cite{xk:4, xkaspm}
 for the necessary condition of discrete decomposability
 of branching laws, 
 where the associated variety of an infinite dimensional representation $\pi$ (an analogue of the characteristic variety of a $\mathcal{D}$-module) plays an important role.
 The references \cite{xkicm, xkshunki} discuss some applications of discrete branching laws.

Loosely, 
Theorem \ref{thm:deco} says that
 if $C_K(K')$ and $\operatorname{AS}_K(\pi)$
 are not \lq{large}\rq\
 then the restriction $\pi|_{G'}$ is discretely decomposable.  
We note that $C_K(K') =\{0\}$ if $K'=K$,
 and consequently, 
 the assumption \eqref{eqn:criterion} is automatically satisfied.  
In this case, 
 Theorem \ref{thm:deco} is nothing but Harish-Chandra's
 admissibility theorem
 ({\cite{xhc}}).  
For any minimal representation $\pi$
 of a reductive group $G$, 
 we get from \cite{xV}
 that $\operatorname{AS}_K(\pi)$
 is one-dimensional, 
 i.e. $\operatorname{AS}_K(\pi)={\mathbb{R}} v$
 or ${\mathbb{R}}_+ v$ where $v$ is the highest root.  
Thus we can expect that there is a rich family of subgroups $G'$ of $G$ for which the restriction of the minimal representation of $G$ decomposes discretely.

\subsection{Space forms of indefinite Riemannian manifolds}
\label{subsec:3.2}

Before applying Theorem \ref{thm:deco}
 to actual branching problems, 
 we review quickly known results
 about the geometry and global analysis 
 on space forms of indefinite-Riemannian manifolds
 (referred also to as pseudo-hyperbolic spaces,
 generalized hyperboloids,
 etc.).  

We set
\begin{alignat*}{2}
     \indefriem +{p,q}:=&\{(x,y) \in {\Bbb R}^{p+1} \oplus {\mathbb R}^q:
                  ||x||^2-||y||^2=1 \}
            &&\simeq O(p+1,q)/O(p,q), 
\\  
     \indefriem -{p,q}:=&\{(x,y) \in {\Bbb R}^{p} \oplus {\mathbb R}^{q+1}:
                  ||x||^2-||y||^2=-1 \}
            &&\simeq O(p,q+1)/O(p,q).  
\end{alignat*}
We note that $X_+^{p,0} \simeq S^p$ and $X_-^{0,q} \simeq S^q$.
By switching the factor,
 we have $\indefriem {+}{p,q} \simeq \indefriem {-}{q,p}$.

We induce an indefinite Riemannian structure on 
 $\indefriem{+}{p,q}$ and $\indefriem{-}{p,q}$
 from the ambient space ${\mathbb{R}}^{p+1,q}$
 and ${\mathbb{R}}^{p,q+1}$, 
 respectively.  
Then, 
$\indefriem{+}{p,q}$
 and $\indefriem{-}{p,q}$
 have constant sectional curvatures.
 Here is a summary of
 indefinite-Riemannian manifolds $X_+^{p,q}$ and $X_-^{p,q}$:
\[
\begin{array}{c|cc}
    & \text{ sectional curvature $\kappa$ } & \text{ signature of metric tensor }
\\[1ex]
\hline
\\
\text{ $\indefriem +{p,q}$ } & \text{ $\kappa \equiv +1$ } & \text{$(p,q)$}
\\[1ex]
\text{ $\indefriem -{p,q}$ } & \text{ $\kappa \equiv -1$ } & \text{ $(p,q)$ }
\end{array}
\]

Let $L^2(\indefriem + {p-1,q})$ be the Hilbert space
 of square integrable functions on $\indefriem + {p-1,q}$
 with respect to the induced volume element.  
For $\lambda \in \Bbb C$,
 we set
$$
     V_{\lambda}^{p,q}
   :=\{f \in L^2(\indefriem +{p-1,q}):
         \widetilde{\Delta}_{\indefriem + {p-1,q}}f=(\frac 1 4-\lambda^2)f \}, 
$$
where the Yamabe operator $\widetilde {\Delta}_{\indefriem + {p,q}}$
 takes the following form:
\begin{equation}
\label{eqn:yamabe}
  \widetilde{\Delta}_{\indefriem {+}{p,q}}=\Delta_{\indefriem{+}{p,q}}-\frac 1 4(p+q)(p+q-2).  
\end{equation}
Clearly, the isometry group $\operatorname{Isom}(\indefriem + {p-1,q}) \simeq O(p,q)$
 preserves $V_{\lambda}^{p,q}$ for any $\lambda \in {\mathbb{C}}$.  
The representations on $V_{\lambda}^{p,q}$ are called 
 {\it{discrete series representations}}
 for $\indefriem+{p-1,q}$ if $V_\lambda^{p,q} \neq \{0\}$, 
 which were studied by Gelfand, Graev, Vilenkin, Shintani, Molchanov, 
 Faraut,  and Strichartz
 among others.  
We summarise:

\begin{prop}
\label{prop:discrep}
\begin{enumerate}
\renewcommand{\labelenumi}{\rm{\arabic{enumi})}}
\item
$(p =1)$\enspace
$V_{\lambda}^{p,q} = \{0\}$
 for any $\lambda \in \Bbb C$.  
\item
$(p \ne 1)$\enspace
$V_{\lambda}^{p,q} \ne \{0\}$
  $\Leftrightarrow$  
 $\lambda \in \frac {p+q}{2} + 2 \Bbb Z$
 and $\lambda \ne 0$.  
\par
Furthermore,
 $O(p,q)$ acts irreducibly on each $V_{\lambda}^{p,q}$,
 when it is non-zero.  
\end{enumerate}
\end{prop}

The resulting representation in Proposition \ref{prop:discrep} 2)
 will be denoted by $\pi_{\lambda}^{p,q}$.  
Since $V_{\lambda}^{p,q}=V_{-\lambda}^{p,q}$,
 we may and do assume $\operatorname{Re} \lambda \ge 0$
 without loss of generality.
By the coherent continuation of $\pi_{\lambda}^{p,q}$
 for $\lambda >0$
 such that $\lambda \in \frac{p+q}{2} + 2 {\mathbb{Z}}$, 
 we can define irreducible unitary representations $\pi_{0}^{p,q}$
 ($p+q$:even)
 and $\pi_{-\frac 1 2}^{p,q}$  ($p+q$:odd) of $O(p,q)$.  
These representations do not lie in $L^2(\indefriem + {p-1,q})$
 but enjoy analogous algebraic properties
 to $\pi_{\lambda}^{p,q}$
 ($\lambda >0$)
 (see \cite[\S 6]{xk:1} or 
 \cite[\S 5.4]{xkors2} 
 the vanishing results on cohomologies in details).

\subsection{Quantative branching laws}
\label{subsec:3.3}

We return to the setting of Section \ref{subsec:2.2}.  
The flat indefinite-Riemannian manifold 
${\mathbb{R}}^{p,q}$
 may be seen as the direct product
 of two flat spaces:
$$
({\mathbb{R}}^p,dx_1^2+ \cdots +dx_p^2)
\times
({\mathbb{R}}^q,-dx_{p+1}^2- \cdots -dx_{p+q}^2).  
$$
Likewise, 
 the direct product of two space forms:
$$
   Y:= \indefriem{+}{p',q'} \times \indefriem{-}{p'',q''}
$$
is locally conformal to ${\mathbb{R}}^{p,q}$
 for any $p',q',p'',q''$
 such that 
$$
 p'+p''=p,
 \qquad
  q'+q''=q.  
$$

This local conformal map is given as follows:
For $u=((\xi_0,\xi',\eta'), (\xi'',\eta'',\eta_0))
\in{\mathbb{R}}^{1+p'+q'} \oplus {\mathbb{R}}^{p''+q''+1}, 
$
 we set
$$
    \Phi(u) := \frac{2}{\xi_0 + \eta_0} (\xi', \eta', \xi'', \eta'').  
$$
Then, 
 the restriction of $\Phi$ to $Y$ is conformal
 (see \cite[Lemma 3.3]{xkors1}, 
 for example).  
More precisely, the map
\begin{equation}
\label{eqn:YR}
     \Phi:\indefriem{+}{p',q'} \times \indefriem{-}{p'', q''} \to {\mathbb{R}}^{p,q}
\end{equation}
 is well-defined and conformal in the open dense set $Y'$ of $Y$, 
 defined by 
 $\xi_0 + \eta_0 \ne 0$.    
Correspondingly, 
 if we set
\begin{equation}
\label{eqn:Phitilde}
     (\widetilde{\Phi}^{\ast} f)(u):= (\frac{2}{\xi_0+ \eta_0})^{\frac{p+q-2}{2}}f(\Phi(u))
\end{equation}
then $\widetilde{\Phi}^{\ast} f$ solves
 $\widetilde{\Delta}_{Y} \widetilde{\Phi}^{\ast} f=0$
 on $Y'$ if $\square_{p,q}f=0$
 (see \cite[Proposition 2.6]{xkors1}).  
Here, 
 $\widetilde {\Delta}_Y$ is the Yamabe operator on $Y$, 
 which amounts to 
\begin{align*}
   {\widetilde{\Delta}}_{Y}
  = &
   {\widetilde{\Delta}}_{\indefriem +{p',q'}}
  - 
   {\widetilde{\Delta}}_{\indefriem -{p'',q''}}
\\
  =& 
   \Delta_{\indefriem +{p',q'}}
   - 
  \Delta_{\indefriem -{p'',q''}} 
  -\frac{1}{4}
   (p'+ q'- p'' -q'')(p+q-2).  
\end{align*}
Hence,
 we can realize the minimal representation $\varpi$ of $O(p+1,q+1)$
on the solution space
$
    \widetilde{\Delta}_{Y}F=0  
$
 as well through $\widetilde{\Phi}^*$.

We note 
 that the map \eqref{eqn:YR} is two to one at generic points.  
In order to give a global action
 of the group $O(p+1, q+1)$
 on the solution space to $\widetilde {\Delta}_Y F =0$, 
 we need to 
define $F = \widetilde{\Phi}^* f$ by \eqref{eqn:Phitilde} for $\xi_0+\eta_0>0$,
and by the parity condition for $\xi_0+\eta_0<0$
 so that $F(-u)=(-1)^{\frac{p-q}{2}}F(u)$ holds 
 (see \cite[(4.4.2a)]{xkors3}).

In light that the isometry group of $Y = \indefriem{+}{p',q'}  \times \indefriem{-}{p'',q''}$
 is the reductive group
 $O(p'+1, q') \times O(p'', q''+1)$,
 it is natural to consider the branching law
 of the minimal representation $\varpi$
 with respect to the following symmetric pair
$$
     O(p+1,q+1) \downarrow O(p'+1,q') \times O(p'',q''+1)
$$
by using the geometric model $Y$.

In this setting, 
 the criterion \eqref{eqn:criterion} of Theorem \ref{thm:deco} holds
 if and only if $p''=0$
 or $q'=0$
 (see \cite[Theorem 4.2]{xkors2}).  
Then, 
it follows from Theorem \ref{thm:deco}
 that $\varpi$ decompose discretely.
For the description of the irreducible decomposition, 
 we define the space of spherical harmonics of degree $l$ by
\begin{align}
\mathcal{H}^l(\mathbb{R}^m) 
:= &\{ \varphi \in C^\infty (S^{m-1}): \Delta_{S^{m-1}}\varphi = -l(l+m-2)\varphi \} \notag \\
= &\left\{ \varphi \in C^\infty (S^{m-1}): \widetilde{\Delta}_{S^{m-1}}\varphi = \left( \frac14 - \left( l + \frac{m-2}{2} \right)^2 \right) \varphi \right\}. \label{eqn:sph}
\end{align}
The orthogonal group $O(m)$ acts irreducibly on $\mathcal{H}^l(\mathbb{R}^m)$ for any $l \in \mathbb{N}$.

Here is the branching law together with quantative information on the invariant inner product:

\begin{thm}[{\rm{(see {\cite[Theorem B]{xkors2}})}}]
\label{thm:twistpull}
Suppose $p+q$ $(>2)$ is even,
 $q=q'+q''$,
 and $p,q >0$.  
Then the twisted pull-back $\widetilde{\Phi}^*$
 of the conformal map $\Phi:Y \to {\mathbb{R}}^{p,q}$
 induces the following quantative branching law:
\begin{enumerate}
\renewcommand{\labelenumi}{\rm{\arabic{enumi})}}
\item
 (branching law; $O(p+1,q+1) \downarrow O(p+1,q') \times O(q'' +1)$).
\begin{equation}
\label{eqn:discdecompopq}
  \varpi^{p+1,q+1}|_{O(p+1,q')\times O(q''+1)}
  \simeq
  {\sum_{l=0}^{\infty}}^{\oplus} \pi_{l+\frac {q''}{2}-\frac 1 2}^{p+1,q'}
  \otimes {\cal H}^l({\Bbb R}^{q''+1})
\end{equation}
Here, 
 the right-hand side of \eqref{eqn:discdecompopq}
 is a multiplicity-free Hilbert direct sum
 of irreducible representations
 of $O(p+1, q') \times O(q'')$.  
\item 
(Parseval-type theorem).
For $f \in \Sol_0(\square_{p,q})$,
 we expand 
 $\widetilde{\Phi}^* f$ into the series $\sum_l F_l$
 according to the discrete decomposition \eqref{eqn:discdecompopq}.  
Then we have:
\begin{equation}
\label{eqn:ppthm}
   ||f||_{{\mathbb{R}}^{p,q}}^2
  =\sum_{l=0}^{\infty} (l+\frac{q''}{2}-\frac 1 2)||F_l||_{L^2(Y)}^2
\end{equation}
Here $||\,\,||_{{\mathbb{R}}^{p,q}}$ is the norm  
 defined in Theorem \ref{thm:definpro}.  
\end{enumerate}
\end{thm}

In view of \eqref{eqn:sph}, 
 the self-adjoint operator $\frac14-\widetilde{\Delta}_{S^{q''}}$ is non-negative,
 and therefore we can define a pseudo-differential operator 
$
     \left(\frac14 - \widetilde{\Delta}_{S^{q''}} \right)^{\frac14}
$ 
  on $Y = X_+^{p,q'} \times S^{q''}$ as well as on $S^{q''}$.

Hence, 
 we get another expression on the invariant inner product
 of the minimal representation in the geometric model $Y$
 by means of a pseudo-differential operator:

\begin{crl}
Suppose $p+q\,(>2)$ is even, $q=q'+q''$, and $p,q''>0$.
Let $\left(\frac14 - \widetilde{\Delta}_{S^{q''}} \right)^{\frac14}$ be the pseudo-differential operator on $Y = X_+^{p,q'} \times S^{q''}$.
We set $F = \widetilde{\Phi}^* f$ for $f \in \overline{\Sol_0(\square_{p,q})}$.
Then
\begin{equation}\label{eqn:pso}
\| f \|_{\mathbb{R}^{p,q}}^2 = 
\left\| \left( \frac14 - \widetilde{\Delta}_{S^{q''}} \right)^{\frac14} F \right\|_{L^2(Y)}^2 .
\end{equation}
\end{crl}

\begin{remark}
\begin{enumerate}
\renewcommand{\labelenumi}{\rm{\arabic{enumi})}}
\item
For $q''=0$ or 1, 
 $l+\frac{q''}{2}-\frac 1 2 \le 0$
 if $l=0$. 
In this case $V_{l+\frac {q''}2 - \frac 1 2}^{p+1, q'}=\{0\}$. 
Nevertheless, 
we can justify the summand in \eqref{eqn:ppthm} by using the argument of the analytic continuation.    

\item
In the case $p''=q'=0$,
 the branching law 
 \eqref{eqn:discdecompopq} is nothing but the $K$-type formula
 of the minimal representation $\varpi$, 
 and \eqref{eqn:ppthm} was proved earlier 
 by Kostant \cite{xKo} for $p=q=3$, 
 and by Binegar and Zierau \cite{xBiZi} for general $p$, $q$
 such that $p+q$ is even and greater than 2.  

\item
In the case $q''=0$,
 we have $Y \simeq \indefriem + {p,q} \times S^0$,
 namely,
 $Y$ consists of two copies of $\indefriem + {p,q}$.  
Then Theorem \ref{thm:twistpull} 
asserts
 that the minimal representation splits into two components, namely,
$$
     \varpi^{p+1,q+1}|_{O(p+1,q)}
     \simeq 
     \pi_{-\frac{1}{2}}^{p+1,q}\oplus\pi_{\frac{1}{2}}^{p+1,q}
$$
 because ${\cal H}^l({\Bbb R}^1)=0$
 for $l \ge 2$.  
\item
In the case $p''=0$ and $p' = q'=1$, 
 we are dealing with the branching law
 for the pair
$$
   O(2,q+1) \downarrow O(2,1) \times O(q).  
$$
The branching law \eqref{eqn:discdecompopq} in this special case
 yields a setting of the deformation
 of the Fourier transform 
 (see Section \ref{sec:4}).  
\end{enumerate}
\end{remark}

\section{Deformation of Fourier transforms}
\label{sec:4}

Minimal representations give us also a hint to define a
generalization of the Fourier transform.
In this section,
we introduce a holomorphic semigroup $\mathcal{I}_{k,a}(z)$ consisting
of Hilbert--Schmidt operators with three parameters:
\begin{itemize}
\item[$a$:]
interpolating minimal representations of simple groups of type $C$ and
 $D$,
\item[$k$:]
Dunkl deformation parameter (multiplicities on the root system),
\item[$z$:]
complex number,
\end{itemize}
such that the \textit{operator-valued boundary value}
\[
\lim_{\operatorname{Re}z\downarrow0} \mathcal{I}_{k,a}(z)
\]
of Hilbert--Schmidt operators 
 yields a one-parameter group of unitary operators.
The underlying idea may be seen as a descendant
 of Sato's hyperfunction theory \cite{xsato}
 and also that of the Gelfand--Gindikin program
\cite{xGeGi,xOls,xStanton}
 for unitary representations
 of real reductive groups.  
We shall see in Diagram \ref{fig:11} 
that the Euclidean Fourier transform,
the Hankel-type transform,
and the Dunkl transform, etc.\ arise naturally as the special values
of 
$
   \mathcal{I}_{k,a}(\frac{\pi i}{2})= \lim_{\varepsilon \downarrow 0}\mathcal{I}_{k,a}(\frac {\pi i}{2}+\varepsilon).
$

\subsection{$L^2$-model of minimal representations}
\label{subsec:4.1}

We return to the setting of Section \ref{subsec:2.2}.  
If a tempered distribution 
$f \in \mathcal{S}'(\mathbb{R}^{p+q})$ 
satisfies the differential equation $\square_{p,q} f = 0$,
then it is easy to see that 
its Fourier transform $\mathcal{F}f$ is supported on
the characteristic variety
\begin{equation}\label{eqn:Xi}
\Xi :=
\{ \xi \in \mathbb{R}^{p+q}:
   \xi_1^2 +\dots+ \xi_p^2 - \xi_{p+1}^2 -\dots- \xi_{p+q}^2 = 0 \}.
\end{equation}
Much more than this,
the following theorem holds
 (see \cite[Theorem 6.2]{xkors2}):
\begin{thm}\label{thm:4.1}
For $p+q>2$ even and $p,q>0$,
the Euclidean
 Fourier transform $\mathcal{F} \equiv \mathcal{F}_{\mathbb{R}^{p+q}}$
induces the bijection: 
\[
\mathcal{F}: \overline{\mathcal{S}ol_0(\square_{p,q})}
\overset{\sim}{\to} L^2(\Xi).
\]
It is an isometry up to the scalar multiplication by 
 $2^{\frac{p+q+2}{2}}\pi^{\frac{p+q+1}{2}}$.
\end{thm}

Here, 
$\overline{\mathcal{S}ol_0(\square_{p,q})}$
is the Hilbert space
with respect to the \textit{conservative quantity}
$(\ ,\ )$ defined in Theorem \ref{thm:definpro},
and $L^2(\Xi)$ denotes the Hilbert space consisting of square
integrable functions with respect to the canonical measure on $\Xi$.
The non-trivial part of Theorem \ref{thm:4.1} is to show
 that $\operatorname{Image} {\mathcal{F}} \cap L^2(\Xi) \ne \{0\}$.
See \cite[Theorem 6.2]{xkors2} for the proof.

It follows from Theorem \ref{thm:4.1} that
 we can realize the minimal representation of the indefinite orthogonal group
$O(p+1,q+1)$ on the Hilbert space $L^2(\Xi)$
(\textit{$L^2$-model})
from the one on $\overline{\mathcal{S}ol_0(\square_{p,q})}$
(\textit{conformal model}).   

At this moment,
 we remark a distinguishing feature of minimal representations 
 (see Appendix in Section \ref{sec:5}).  
Unlike well-understood family of irreducible unitary representations
of real reductive  groups such as unitary principal series
representations or discrete series representations,
minimal representations are too `small' that there is no existing
geometric model for which both group actions and the Hilbert structure
are given in a simple manner
 (cf. \cite{xbk, xTo}).
We pin down the advantages of the aforementioned two models:
\[
\begin{array}{c|cc}
    & \text{ Group action } & \text{ Hilbert structure }
\\[1ex]
\hline
\\
\text{ Conformal model $\overline{\mathcal{S}ol(\square_{p,q})}$ } 
& \text{ simple } & \text{\textcircled{1}}
\\[1ex]
\text{ $L^2$-model $L^2(\Xi)$ } 
& \text{\textcircled{2}} & \text{ simple }
\end{array}
\]

Finding
 the missing parts \textcircled{1}
and \textcircled{2} is
 interesting, 
 particularly because it interacts with other branches of mathematics.  
Representation theoretic consideration plays a guiding
role in formalizing problems there.
In fact, 
 we have seen in Theorem \ref{thm:definpro} 
that \textcircled{1} brought us to the construction of  conservative
quantities for ultra-hyperbolic equations,
whereas \textcircled{2} leads us to the notion of a Fourier
transform on the isotropic cone $\Xi$ \cite{xbkocras,xkmano2,xkmanoAMS},
as discussed below.

{}From now,
we consider the missing part \textcircled{2}.
In order to find the global formula of group actions on the
$L^2$-model, 
let us clarify what is trivial and what will be the crucial operator.
We observe that there is a maximal parabolic subgroup $P$ of 
$G=O(p+1,q+1)$ that
contains the conformal transformation group
\[
\operatorname{Conf}(\mathbb{R}^{p,q})
\simeq
(\mathbb{R}_{>0} \times O(p,q))
\ltimes \mathbb{R}^{p+q}
\]
as a
subgroup of index two.
Then we have the Bruhat decomposition
\[
G = P \amalg P w P,
\]
where
$w = \begin{pmatrix} I_{p+1} & 0 \\ 0 & -I_{q+1} \end{pmatrix}$.
In fact, 
 the Euclidean Fourier transform ${\mathcal {F}}_{{\mathbb{R}}^N}$
 appear as the {\it{unitary inversion operator}}
 of the Segal--Shale--Weil representation 
 of the metaplectic group $Mp(N,{\mathbb{R}})$,
 which is also a minimal representation.  
See \cite[Chapter 1]{xkohcrcras} 
 for the comparison of ${\mathcal{F}}_{\Xi}$ and ${\mathcal{F}}_{{\mathbb{R}}^N}$
 from this point of view.

In the $L^2$-model of the minimal representation $\pi$ 
 of $G$ on
$L^2(\Xi)$,
the $P$-action is simple,
namely, it is given just by translations and multiplications \cite{xkors3}.
Hence, it is enough to
find the single unitary operator
(\textit{unitary inversion operator}) $\pi(w)$ 
in order to fill the missing part \textcircled{2}.
We set
\begin{equation}\label{eqn:FXi}
\mathcal{F}_\Xi := c\pi(w),
\end{equation}
where $c$ is the phase factor.
Algebraically, $\mathcal{F}_\Xi$ intertwines the multiplication of
coordinate functions $\xi_j$ $(1 \le j \le p+q)$ with the Bargmann--Todorov
operators $R_j$ $(1 \le j \le p+q)$ which are mutually commuting 
differential operators of second order on $\Xi$
(see \cite{xbt}, \cite[Chapter 1]{xkmanoAMS}).

This algebraic feature is similar to the classical fact that the Euclidean Fourier transform
$\mathcal{F}_{\mathbb{R}^N}$ intertwines the multiplication operators
$\xi_j$ and the differential operators $\sqrt{-1}\partial_j$
$(1 \le j \le N)$.

The goal of this section is to explain these operators
$\mathcal{F}_\Xi$ and $\mathcal{F}_{\mathbb{R}^N}$ from a broader
point of view,
by constructing continuous family of operators that include  
$\mathcal{F}_\Xi$ and $\mathcal{F}_{\mathbb{R}^N}$ as their special
values.

For this,
we limit ourselves to the case $p=1$.
Then,
the light cone $\Xi$ (see \eqref{eqn:Xi}) splits into the forward light cone
$\Xi_+$ and the backward light cone $\Xi_-$
 according as $x_1>0$ and
$x_1<0$.
The unitary inversion operator $\mathcal{F}_\Xi$ preserves the direct
sum
\begin{equation}
\label{eqn:Xi2}
L^2(\Xi) = L^2(\Xi_+) \oplus L^2(\Xi_-).
\end{equation}
For later purpose,
 we set $q=N$.
Then the projection to the second factor,
$\mathbb{R}^1 \oplus \mathbb{R}^N \to \mathbb{R}^N$,
 induces the following isomorphism between the Hilbert spaces:
\begin{equation}\label{eqn:LXi}
L^2(\Xi_+) \simeq L^2(\mathbb{R}^N, \|x\|^{-1} dx).
\end{equation}
Via \eqref{eqn:LXi},
the unitary inversion operator $\mathcal{F}_\Xi$ 
on $L^2(\Xi_+)$ may be seen as a unitary operator on
$L^2(\mathbb{R}^N, \|x\|^{-1} dx)$.
The explicit formula
 of ${\mathcal{F}}_{\Xi}$ in the coordinates of ${\mathbb{R}}^N$
 was given in \cite{xkmano2}.  
In this framework,
we can construct a holomorphic family of bounded operators so that
the unitary operator ${\mathcal {F}}_{\Xi}$ is obtained
 as the limit of
holomorphic objects.  
Deformation in the Dunkl setting \cite{xbkocras, xbko}
 is also built on this formulation.  
We will discuss those operators in this generality in Section \ref{subsec:4.4}.  

An alternative approach was taken in \cite{xkmanoAMS} based on the
Barnes--Mellin integral to find the kernel function of
$\mathcal{F}_\Xi$ for general $p$ and $q$.

\subsection{Hermite semigroup and Fourier transform}
\label{subsec:4.2}

We begin with recalling a general fact
 on the classical Hermite operator on $\mathbb{R}^N$
 (e.g.\ \cite{xfolland,Howe}):

\begin{equation}\label{eqn:Her}
\Delta - \|x\|^2
= \sum_{j=1}^N \frac{\partial^2}{\partial x_j^2}
  - \sum_{j=1}^N x_j^2.
\end{equation}
Then, 
 $\Delta-\|x\|^2$ extends to a self-adjoint operator on
$L^2(\mathbb{R}^N)$.  
We normalize the Euclidean Fourier transform $\mathcal{F}_{\mathbb{R}^N}$
on $L^2(\mathbb{R}^N)$ as 
\[
(\mathcal{F}_{\mathbb{R}^N}f)(\xi)
=
\frac{1}{(2\pi)^{\frac{N}{2}}}
\int_{\mathbb{R}^N} f(x) e^{-i\langle x,\xi\rangle} dx.  
\]
Then, 
${\mathcal{F}}_{{\mathbb{R}}^N}$ is written as a special value
of the one-parameter group of unitary operators
\[
\chi(t) :=
\exp\left( \frac{it}{2} (\Delta - \|x\|^2)\right), 
\]
namely, 
 we have 
\begin{equation}
\label{eqn:Fexp}
\mathcal{F}_{\mathbb{R}^N}
=  e^{\frac{1}{4}\pi i N}
\exp \left( \frac{\pi i}{4} (\Delta - \|x\|^2) \right).
\end{equation}
Further, 
the one-parameter group $\chi(t)$ of unitary operators
 extends to a holomorphic semigroup
$I(z)$ defined by
\begin{equation}\label{eqn:Hsemi}
I(z) = \exp \frac{z}{2} (\Delta - \|x\|^2)
\quad\text{for $\operatorname{Re} z > 0$}.
\end{equation}

The semigroup $I(z)$ is called the {\it{Hermite semigroup}},
 and it is expressed as an integral transform against the Mehler kernel 
 \cite{xfolland, Howe},
 a Gaussian type kernel.  

Next, we consider another differential operator on ${\mathbb{R}}^N$, 
\[
\|x\| \Delta - \|x\|.  
\]
It turns out
 that this operator has a self-adjoint extension on the Hilbert space
$L^2(\mathbb{R}^N,\|x\|^{-1} dx)$
(see \cite[Section 1.1]{xkmano2}).
Moreover,
 an analogous formula to \eqref{eqn:Fexp} holds:
 via the identification \eqref{eqn:LXi},
the `Fourier transform' $\mathcal{F}_{\Xi}$ on the forward light cone
$\Xi_+$ can be expressed as
\begin{equation}
\label{eqn:Fexp2}
\mathcal{F}_\Xi
= c \exp \left( \frac{\pi i}{2} (\|x\| \Delta - \|x\|) \right),
\end{equation}
where $c = e^{\frac{1}{2}\pi i (N-1)}$ is the phase factor.
Then, 
 the expression \eqref{eqn:Fexp2} allows us to see 
 $\mathcal{F}_\Xi$ as the limit of the following holomorphic semigroup
(\textit{Laguerre semigroup})
\begin{equation}\label{eqn:Lsemi}
\mathcal{I}(z) = \exp z (\|x\| \Delta - \|x\|),
\quad\text{for $\operatorname{Re}z>0$}, 
\end{equation}
as $z \to \frac {\pi i}{2}+0$.  
The kernel function of ${\mathcal{I}}(z)$ 
 is given in terms of the Bessel function \cite{xkmano1}.

Interpolating $\Delta - \|x\|^2$ and 
$\|x\| \Delta - \|x\|$, 
 namely, 
 the infinitesimal generators
 of the Hermite semigroup \eqref{eqn:Hsemi} and the Laguerre
semigroup \eqref{eqn:Lsemi},
 we consider the differential operator
\[
\Delta_{0,a}
:= \|x\|^{2-a} \Delta - \|x\|^a.
\]
It might not be so obvious that the symmetric operator $\Delta_{0,a}$ has a
self-adjoint extension on the Hilbert space
$L^2(\mathbb{R}^N,\|x\|^{a-2}dx)$.  
In fact, 
 it is the case.  
The proof uses representation theory
(see Proposition \ref{thm:F}), 
 and the same idea works in a more general setting  
 of Dunkl's differential-difference operators. 
Thus,
we shall introduce a holomorphic semigroup
$\mathcal{I}_{k,a}(z)$ 
 with infinitesimal generator $\Delta_{k,a}$
 (see \eqref{eqn:kaLap} below for the definition)
 for $\operatorname{Re}z>0$ with parameters $k$
and $a$ in Section \ref{subsec:4.3}.

\begin{figure}[H]
\renewcommand{\figurename}{Diagram}
 \renewcommand{\thefigure}{\thesubsection}
\setlength{\unitlength}{.8mm}
\begin{picture}(160,110)
\put(26,101){\makebox(100,10){\fbox{$(k,a)$-generalized
                              Fourier transform $\mathcal{F}_{k,a}$}}} 
\put(12,80){$\scriptstyle a \to 2$}
\put(5.5,46){%
\setlength{\unitlength}{0.0007in}
\begingroup\makeatletter\ifx\SetFigFont\undefined%
\gdef\SetFigFont#1#2#3#4#5{%
  \reset@font\fontsize{#1}{#2pt}%
  \fontfamily{#3}\fontseries{#4}\fontshape{#5}%
  \selectfont}%
\fi\endgroup%
{\renewcommand{\dashlinestretch}{30}
\begin{picture}(1530,2431)(0,-10)
\put(3172.000,-292.000){\arc{6328.507}{3.2365}{4.1638}}
\path(5.967,130.650)(22.000,8.000)(65.571,123.765)
\end{picture}
}
}
\put(136,80){$\scriptstyle a \to 1$}
\put(111,46){%
\setlength{\unitlength}{0.0007in}
\begingroup\makeatletter\ifx\SetFigFont\undefined%
\gdef\SetFigFont#1#2#3#4#5{%
  \reset@font\fontsize{#1}{#2pt}%
  \fontfamily{#3}\fontseries{#4}\fontshape{#5}%
  \selectfont}%
\fi\endgroup%
{\renewcommand{\dashlinestretch}{30}
\begin{picture}(1529,2431)(0,-10)
\put(-1642.000,-292.000){\arc{6328.507}{5.2609}{6.1882}}
\path(1464.429,123.765)(1508.000,8.000)(1524.033,130.650)
\end{picture}
}
}
\put(75,95){$\bigg\uparrow$}
\put(77,95){$\scriptstyle z \to \frac{\pi i}{2}$}
\put(51,79){\framebox(64,11){{\framebox(61,8){holomorphic semigroup
                             $\mathcal{I}_{k,a}(z)$}}}}
\put(45,71){$\scriptstyle a \to 2$}
\put(44,70){\rotatebox[origin=b]{30}{$\xleftarrow{\kern19mm}$}}
\put(90,70){\rotatebox[origin=b]{-30}{$\xrightarrow{\kern19mm}$}}
\put(103,71){$\scriptstyle a \to 1$}
\put(29,59){\fbox{$\mathcal{I}_{k,2}(z)$}}
\put(108,59){\fbox{$\mathcal{I}_{k,1}(z)$}}
\put(16,51){$\scriptstyle z \to \frac{\pi i}{2}$}
\put(20,50){\rotatebox[origin=b]{30}{$\xleftarrow{\kern10mm}$}}
\put(43,50){\rotatebox[origin=b]{-30}{$\xrightarrow{\kern10mm}$}}
\put(51,51){$\scriptstyle k \to 0$}
\put(97,51){$\scriptstyle k \to 0$}
\put(99,50){\rotatebox[origin=b]{30}{$\xleftarrow{\kern10mm}$}}
\put(122,50){\rotatebox[origin=b]{-30}{$\xrightarrow{\kern10mm}$}}
\put(130,51){$\scriptstyle z \to \frac{\pi i}{2}$}
\put(3,34){\framebox(33,11){\begin{tabular}{c}
                            Dunkl transform\\[-.5ex] $\mathcal{D}_k$ 
                            \end{tabular}}}
\put(38,34){\framebox(41,11){\begin{tabular}{c}
                            Hermite semigroup\\[-.5ex] $I(z)$ 
                            \end{tabular}}}
\put(83,34){\framebox(41,11){\begin{tabular}{c}
                              Laguerre semigroup 
                              \end{tabular}}}
\put(127,34){\framebox(24,11){\begin{tabular}{c}
                            $\mathcal{F}_{k,1}$
                            \end{tabular}}}
\put(18,27){$\scriptstyle k \to 0$}
\put(20,27.5){\rotatebox[origin=b]{-30}{$\xrightarrow{\kern10mm}$}}
\put(43,27.5){\rotatebox[origin=b]{30}{$\xleftarrow{\kern10mm}$}}
\put(51,27){$\scriptstyle z \to \frac{\pi i}{2}$}
\put(94,27){$\scriptstyle z \to \frac{\pi i}{2}$}
\put(99,27.5){\rotatebox[origin=b]{-30}{$\xrightarrow{\kern10mm}$}}
\put(122,27.5){\rotatebox[origin=b]{30}{$\xleftarrow{\kern10mm}$}}
\put(131,27){$\scriptstyle k \to 0$}
\put(20,16){\framebox(37,7){Fourier transform}}
\put(99,16){\framebox(37,7){Hankel transform}}
\put(37,5){\rotatebox{90}{\hbox to 2.5em{\dotfill}}}
\put(48,10){$\Leftarrow$ `unitary inversion operator' $\Rightarrow$}
\put(117,5){\rotatebox{90}{\hbox to 2.5em{\dotfill}}}
\put(5,0){\parbox{50mm}{\begin{tabular}{c}the Weil representation of\\[-.5ex]
                        $Mp(N,\mathbb{R})$\end{tabular}}}
\put(90,0){\parbox{50mm}{\begin{tabular}{c}the minimal representation of\\[-.5ex]
                         $O(2,N+1)$\end{tabular}}}
\end{picture}
\label{fig:11}
\caption{Special values of holomorphic semigroup $\mathcal{I}_{k,a}(z)$}
\end{figure}
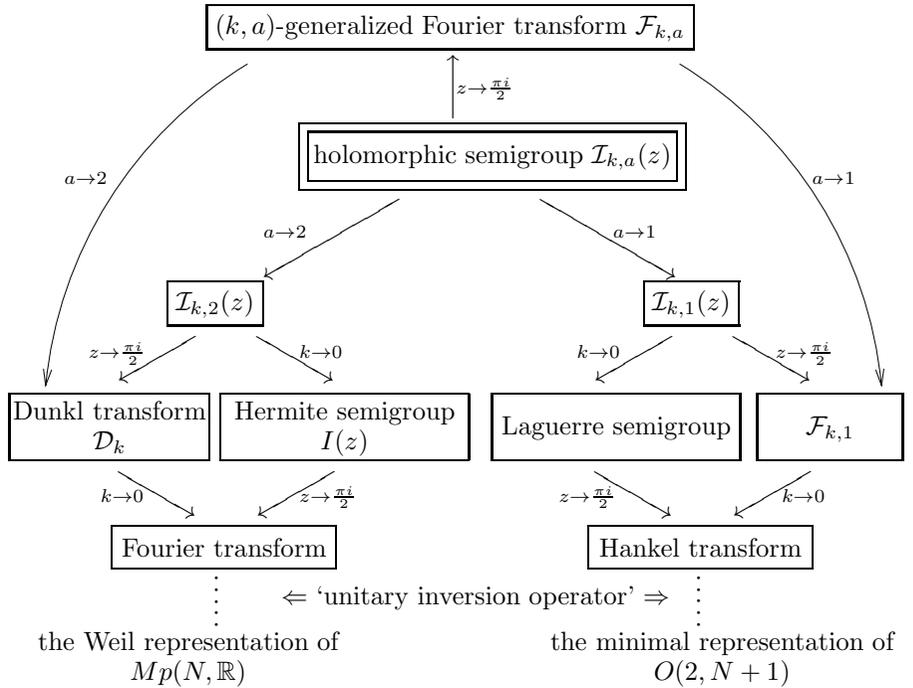

In Diagram \ref{fig:11},
we have summarised some of the deformation properties by indicating
the limit behaviour of the holomorphic semigroup 
$\mathcal{I}_{k,a}(z)$.
The specialization
$\mathcal{I}_{k,a}(\frac{\pi i}{2})$ leads us to a 
$(k,a)$-\textit{generalized Fourier transform}
$\Fka$ (up to a phase factor),
which reduces to the Fourier transform
($a=2$ and $k\equiv0$),
the Dunkl transform $\mathcal{D}_k$
($a=2$ and $k\equiv0$),
and the Hankel-type transform
($a=1$ and $k\equiv0$).

\subsection{Holomorphic semigroup $\mathcal{I}_{k,a}(z)$ with two
parameters $k$ and $a$}
\label{subsec:4.3}

This subsection introduces a holomorphic semigroup, 
 denoted by $\mathcal{I}_{k,a}(z)$, 
 of which the infinitesimal generator is a self-adjoint differential-difference operator.

Let $\mathfrak{C}$ be the Coxeter group associated with a reduced root system
$\mathcal{R}$ in $\mathbb{R}^N$.
For a $\mathfrak{C}$-invariant function $k\equiv(k_\alpha)$
(\textit{multiplicity function}) on $\mathcal{R}$, 
we set
$$
\langle k\rangle := \frac{1}{2} \sum_{\alpha\in\mathcal{R}} k_\alpha,
$$
and 
write $\Delta_k$ for the Dunkl Laplacian
 on $\mathbb{R}^N$ (see \cite{H}).
This is a differential-difference operator, 
 which reduces to the Euclidean Laplacian 
$\Delta$ when $k \equiv 0$.  

We take $a>0$ to be yet another deformation parameter, 
and define
\begin{equation}\label{eqn:kaLap}
\Delta_{k,a}
 := \Vert x \Vert^{2-a} \Delta_k - \Vert x \Vert^a.
\end{equation}
We define a density on ${\mathbb{R}}^N$ by 
$$
\vartheta_{k,a}(x)
 := \|x\|^{a-2} \prod_{\alpha\in\mathcal{R}}
                   | \langle \alpha,x \rangle |^{k_\alpha}.
$$
The volume of the unit ball with respect to the measure
 $\vartheta_{k,a}(x) d x$ is explicitly known in terms of the gamma function
 owing to the work by Selberg, 
 Macdonald, 
 Heckman, and Opdam
 among others
 but we do not go into details
 (see Etingov \cite{xEt}).

In the case $a=2$ and $k\equiv0$,
 $\vartheta_{0,2}(x) \equiv 1$ and $\Delta_{0,2}$ is the
Hermite operator \eqref{eqn:Her} on $\mathbb{R}^N$.

Here are  basic properties of our differential-difference operator
$\Delta_{k,a}$:

\begin{thm}[{\rm{(see {\cite[Theorem A]{xbko}})}}]\label{thm:A}
Assume $a>0$ and $a+2\langle k\rangle+N-2>0$.
\begin{enumerate}[\upshape 1)]
\item  
$\Delta_{k,a}$ extends to a self-adjoint operator on 
the Hilbert space $L^2(\mathbb{R}^N,\vartheta_{k,a}(x)dx)$.
\item  
There is no continuous spectrum of $\Delta_{k,a}$. 
\item 
All the discrete spectrum of $\Delta_{k,a}$ is negative.
\end{enumerate}
\end{thm}

We introduce the following operators on $L^2(\mathbb{R}^N,\vartheta_{k,a}(x)dx)$ by
\begin{equation}\label{eqn:02}
\mathcal{I}_{k,a}(z)
:=
 \exp\Bigl(\frac{z}{a}\Delta_{k,a}\Bigr),
\quad
  \text{ for }
  \operatorname{Re} z \ge 0.  
\end{equation}
Correspondingly to the properties of the infinitesimal generator
$\frac{1}{a}\Delta_{k,a}$ in Theorem \ref{thm:A},
we get:
\begin{thm}
[{\rm{(see {\cite[Theorem B]{xbko}})}}]
\label{thm:B}
Retain the assumption of Theorem \ref{thm:A}.
\begin{enumerate}[\upshape 1)]
\item  
$\{\mathcal{I}_{k,a}(z): \operatorname{Re}z>0\}$ forms
a holomorphic semigroup 
in the complex right-half plane $\{z \in {\mathbb{C}}: \operatorname{Re} z >0\}$
 in the sense that ${\mathcal {I}}_{k,a}(z)$ is a Hilbert--Schmidt operator
 on $L^2({\mathbb{R}}^N, {\vartheta}_{k,a}(x) dx)$
 satisfying 
$$
     {\mathcal{I}}_{k,a}(z_{1}) \circ {\mathcal{I}}_{k,a}(z_{2})= {\mathcal{I}}_{k,a}(z_{1}+ z_{2}), 
    \quad
   (\operatorname{Re} z_1, \operatorname{Re} z_2 >0), 
$$
and that the scalar product $({\mathcal{I}}_{k,a}(z) f, g)$
 is a holomorphic function of $z$ for $\operatorname{Re} z>0$,
 for any $f,g \in L^2({\mathbb{R}}^N, {\vartheta}_{k,a} (x) d x)$.   
\item  
$\mathcal{I}_{k,a}(z)$ is a one-parameter group of unitary operators
on the imaginary axis $\operatorname{Re}z=0$.
\end{enumerate}
\end{thm}
We shall call $\mathcal{I}_{k,a}(z)$ as the
 \textit{$(k,a)$-generalized Laguerre semigroup} $\mathcal{I}_{k,a}(z)$.
We note that 
$\mathcal{I}_{0,2}(z)$ is the Hermite semigroup \eqref{eqn:Hsemi} (see \cite{xfolland,Howe}), 
 and ${\mathcal{I}}_{0,1}(z)$ is the Laguerre semigroup
 \eqref{eqn:Lsemi} (see \cite{xkmano1}).

\subsection{$(k,a)$-generalized Fourier transforms
$\mathcal{F}_{k,a}$}
\label{subsec:4.4}

Theorem \ref{thm:B} 2) asserts that
the `boundary value' of the holomorphic
semigroup $\mathcal{I}_{k,a}(z)$
 produces a one-parameter family of
unitary operators. 

The case $z=0$ gives the identity operator, namely, 
$\mathcal{I}_{k,a}(0) = \operatorname{id}$.
The particularly interesting case is when $z = \frac{\pi i}{2}$.
We set 
$$
c := \exp (i \pi \frac{N + 2\langle k \rangle +a-2}{2a})
\quad\text{(phase factor)}, 
$$ 
and define the 
\textit{$(k,a)$-generalized Fourier transform} by
$$
\index{Fka@$\Fka$}%
\Fka
 := c \, \mathcal{I}_{k,a} \Bigl(\frac{\pi i}{2}\Bigr)
= c \exp 
\Bigl( \frac{\pi i}{2a} (\Vert x\Vert^{2-a} \Delta_k - \Vert x\Vert^a) \Bigr).
$$
Then, this operator $\Fka$ for general $a$ and $k$ 
satisfies the following significant
properties:
\begin{thm}[{\rm{({\cite[Theorem D]{xbko}})}}]
\label{thm:D}
Retain the assumption 
 of Theorem \ref{thm:A}.  
\begin{enumerate}[\upshape 1)]
\item
$\Fka$ is a unitary operator on
$L^2(\mathbb{R}^N,\vartheta_{k,a}(x)dx)$.
\item
$\Fka \circ H_{k,a}
 = - H_{k,a} \circ \Fka$.
\enspace
Here, $H_{k,a}$ is the differential operator of first order defined in
\eqref{eqn:Ha}. 
\item
$\Fka \circ \Vert x\Vert^a
= - \Vert x\Vert^{2-a} \Delta_k \circ \Fka$,

$\Fka \circ (\Vert x\Vert^{2-a} \Delta_k)
= -\Vert x\Vert^a \circ \Fka$.
\item
$\Fka$ is of finite order if and only if $a\in\mathbb{Q}$.
Its order is $2m$ if $a=\frac{m}{n}$ such that $(m,n)=1$.
In particular,
$\mathcal{F}_{k,1}$ is of order $2$,
and $\mathcal{F}_{k,2}$ is of order $4$.
\end{enumerate}
\end{thm}

As indicated in Diagram \ref{fig:11},
$\Fka$ reduces to the Euclidean Fourier transform $\mathcal{F}$  on $\mathbb{R}^N$
 if $k\equiv0$ and $a=2$;
to the Dunkl transform $\mathcal{D}_k$ introduced by C. Dunkl himself
if $k>0$ and $a=2$.
The unitary operator 
$
\mathcal{F}_{0,1}
$
arises as the unitary inversion operator 
${\mathcal{F}}_{\Xi}$ on $L^2(\Xi_+)$
 of the minimal representation of the conformal group 
 (see Section \ref{subsec:4.1}).

Our study also contributes to the theory of special functions,
in particular orthogonal polynomials;
indeed we derive several new identities,
for example,
the $(k,a)$-deformation of the classical Bochner--Hecke identity where
the Gaussian function and harmonic polynomials in the classical
setting
($k\equiv0$ and $a=2$)
are replaced respectively with
$\exp(-\frac{1}{a} \|x\|^a)$
and polynomials annihilated by the Dunkl Laplacian.
The $(k,a)$-generalized
Fourier transform $\Fka$ also satisfies a Heisenberg-type inequality.  
This generalizes the classical case
($k\equiv0$ and $a=2$)
and R\"{o}sler's Heisenberg inequality \cite{xros}
($k>0$ and $a=2$).
We refer to \cite{xbko} for full details.

\subsection{Hidden symmetries in the Hilbert space $L^2(\mathbb{R}^N,\vartheta_{k,a}(x)dx)$}
\label{subsec:4.5}

The key idea of the proof for Theorem \ref{thm:4.1}, 
 \ref{thm:A}, and \ref{thm:B}
 is to use more operators rather than the single operator
$\Delta_{k,a}$, 
 and then to appeal representation theory
  of ${\mathfrak {sl}}_2$,
 in particular, 
 the theory of discretely decomposable
 unitary representations.   

\begin{lmm}
\label{lem:E}
Let $k$ be a multiplicity-function on a root system,
and $a \in \mathbb{C}^\times$.
Then, 
the following differential-difference operators on $\mathbb{R}^N \setminus \{0\}$
\begin{align}\label{eqn:Ha}
&E_{k,a}^+:={i\over a} \Vert x\Vert^a,
\nonumber
\\
&E_{k,a}^-:={i\over a}\Vert
x\Vert^{2-a} \Delta_k,
\nonumber
\\
& H_{k,a}:={2\over a} \sum_{i=1}^N x_i
\partial_i+{{N+2\langle k\rangle+a-2}\over a}
\end{align}
form an
$\mathfrak{sl}_2$-triple, namely,
we have:
\[
[H_{k,a}, E_{k,a}^+] = 2E_{k,a}^+,
\quad
[H_{k,a}, E_{k,a}^-] = -2E_{k,a}^-,
\quad
[E_{k,a}^+, E_{k,a}^-] = H_{k,a}.
\]
\end{lmm}

Special cases of Lemma \ref{lem:E} was previously known:
the case $k\equiv 0$ and $a=2$ is the classical harmonic
$\mathfrak{sl}_2$-triple (e.g.\ Howe \cite{Howe}),
the case $k>0$ and $a=2$ by Heckman \cite{H},
and $k\equiv 0$ and  $a=1$ by Kobayashi and Mano \cite{xkmano1}.
The operator 
$\Delta_{k,a}$ (see \eqref{eqn:kaLap}) takes the form,
$$
\Delta_{k,a} = ai (E_{k,a}^+ - E_{k,a}^-),
$$
which may be seen as an element of $\mathfrak{sl}(2,\mathbb{C})$.

Lemma \ref{lem:E} fits well into the framework of discretely decomposable 
representations of reductive groups 
\cite{xk:2,xk:3,xk:4}
 as we discussed in Section \ref{subsec:3.1}: 
\begin{prop}[{\rm{(see {\cite[Theorem 3.31]{xbko}})}}]
\label{thm:F}
If $a>0$ and $a+2\langle k\rangle+N-2>0$,
then the representation of\/ $\mathfrak{sl}(2,\mathbb{R})$ lifts to a unitary representation of
the simply-connected group on $L^2(\mathbb{R}^N,\vartheta_{k,a}(x)dx)$.
The resulting unitary representation is discretely decomposable,
and commutes the obvious
action of the Coxeter group $\mathfrak{C}$.
\end{prop}

This unitary representation plays the central role in the key
 to the proof of Theorems \ref{thm:A}, \ref{thm:B} and \ref{thm:D}. 
An explicit formula of the irreducible decomposition of
$L^2(\mathbb{R}^N,\vartheta_{k,a}(x))$
is found in \cite[Theorem 3.28]{xbko}.
In the special cases $k\equiv0$ and $a=1$ or $2$,
this formula may be regarded as the branching law of the minimal
representations of
$O(2,N+1)^{\widetilde{}}$ \
or $Mp(N,\mathbb{R})$,
respectively (see Diagram \ref{fig:Hs} below).
Correspondingly,
 all the spectrum of $\Delta_{k,a}$ is also obtained explicitly.

In the case $a=2$, 
 the ${\mathfrak {sl}}_2$-action also induces the representation of $SL(2,{\mathbb{C}})$
 on the algebra generated by Dunkl's operators,
 multiplication operators,
 and the Coxeter group.  
The restriction of this action to $SL(2,{\mathbb{Z}})$
 coincides with a special case of the $SL(2,{\mathbb{Z}})$-action
 discovered by Cherednik \cite{xcherednik}
 on the (degenerate) rational DAHA (double affine Hecke algebra).

Theorem \ref{thm:F} asserts that the Hilbert
space
$L^2(\mathbb{R}^N,\vartheta_{k,a}(x)dx)$ has a symmetry of the direct product
group $\mathfrak{C} \times \widetilde{SL(2,{\mathbb{R}})}$ for all $k$ and $a$.
This symmetry becomes larger for special values of $k$ and $a$ as
 below:
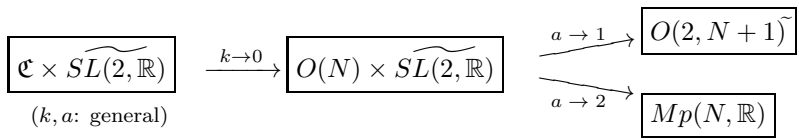
\begin{figure}[H]
\renewcommand{\figurename}{Diagram}
 \renewcommand{\thefigure}{\thesubsection}
\begin{center}
\setlength{\unitlength}{1mm}
\begin{picture}(100,19)(0,2)
\put(84,15){\fbox{$O(2,N+1)^{\widetilde{}}$}}
\put(72,15){\scriptsize $a\to1$}
\put(70,13){\rotatebox[origin=b]{10}{$\xrightarrow{\kern11mm}$}}
\put(0,10){\fbox{$\mathfrak{C}\times\widetilde{SL(2,\mathbb{R})}$}}
\put(26,10){$\xrightarrow{\ k\to 0\ }$}
\put(37,10){\fbox{$O(N)\times\widetilde{SL(2,\mathbb{R})}$}}
\put(3,4){\footnotesize ($k,a$: general)}
\put(70,8){\rotatebox[origin=b]{-10}{$\xrightarrow{\kern11mm}$}}
\put(84,4){\fbox{$Mp(N,\mathbb{R})$}}
\put(72,6){\scriptsize $a\to2$}
\end{picture}
\end{center}
\label{fig:Hs}
\caption{Hidden symmetries in $L^2(\mathbb{R}^N,\vartheta_{k,a}(x)dx)$}
\end{figure}

For $a=2$, this symmetry is given by the Segal--Shale--Weil
representation of the metaplectic group
$Mp(N,\mathbb{R})$.
For $a=1$,
 it is given by the irreducible unitary representation of
the double covering
$O(2,N+1)^{\widetilde{}}$ \ 
of the conformal group
on
$L^2(\mathbb{R}^N,\Vert x\Vert^{-1}dx)$.
Here, 
 as we saw in Theorem \ref{thm:4.1}, 
 we do not need to take a double cover when $N$ is odd.  
Both of them are minimal representations
 and, in particular, they attain the minimum of their Gelfand--Kirillov
dimensions
 among the unitary dual.
In this sense, our continuous parameter $a>0$ interpolates the
$L^2$-models of two minimal
representations of different reductive groups by 
keeping smaller
symmetries 
$O(N) \times \widetilde{SL(2,\mathbb{R})}$.
In view of Lemma \ref{lem:E}, 
 the $(k,a)$-generalized Fourier transform $\Fka$  
($k \equiv 0$, $a=1,2$) arise as
the unitary operators (up to phase factors) corresponding to the
 following element,
$$
  \exp \frac {\pi}{2}\begin{pmatrix} 0 & 1 \\ -1 & 0\end{pmatrix} \in \widetilde {SL(2,{\mathbb{R}})}.  
$$

\section{Appendix --- representation theoretic properties of $\varpi$}
\label{sec:5}
For the reader's convenience,
 we list some representation theoretic properties
 of the irreducible unitary representation $\varpi$
 of the indefinite orthogonal group $G=O(p+1,q+1)$, 
 on which geometric analysis is the motif throughout this article.  

In what follows,
 we assume
$$
     p,q \ge 1
    \,\,
    \text{ and }
    \,\,
    p+q
    \text{ is an even integer greater than two.}
$$
We write $K \simeq O(p+1) \times O(q+1)$
 for a maximal compact group of $G$, 
 ${\mathfrak {g}} \simeq {\mathfrak {o}}(p+1,q+1)$
 for the Lie algebra of $G$,
 ${\mathfrak {g}}_{\mathbb{C}}\simeq {\mathfrak {o}}(p+q+2, {\mathbb{C}})$
 for its complexification,
 and 
$
     {\mathfrak{g}}_{\mathbb{C}}={\mathfrak{k}}_{\mathbb{C}} + {\mathfrak{p}}_{\mathbb{C}}
$ 
 for the complexified Cartan decomposition.  

Here are some properties of $\varpi$ from representation theoretic viewpoints.

\begin{enumerate}
\item[1)]
$\varpi$ is an irreducible unitary representation of $G$.  

\item[2)]
(minimal representation)
\enspace
The representation $\varpi$ is a minimal representation
 in the sense that the annihilator of the underlying $({\mathfrak {g}}_{\mathbb{C}}, K)$-module
 $\varpi_K$
 in the universal enveloping algebra
 $U({\mathfrak {g}}_{\mathbb{C}})$
 is the Joseph ideal if $p+q \ge 6$
 (\cite{xBiZi, xKo}).  
See \cite{xgansavin} for algebraic aspects
 of minimal representations
 of reductive groups
 and the definition
 of the Joseph ideal.  

\item[3)]
(restriction to the identity component)
\enspace
The group $G$ has four connected components,
 and we write $G_0=SO_0(p+1,q+1)$
 for the identity component.  
Then, 
$\varpi$ stays irreducible 
 when restricted to $G_0$
 if and only if $p,q>1$.  
\item[4)]
(highest weight module case)\enspace
If $p=1$ or $q=1$, 
 then the restriction $\varpi|_{G_0}$ is a direct sum of two irreducible representations, 
 one is a highest weight representation $\varpi_+$ 
 and the other is a lowest weight representation $\varpi_-$. 
As we have seen in \eqref{eqn:Xi2}, 
 this decomposition $\varpi|_{G_0}= \varpi_+ \oplus \varpi_-$
 corresponds to the direct sum
$$
     L^2(\Xi) =L^2(\Xi_+) \oplus L^2(\Xi_-)
$$
in the $L^2$-model.  
Both $\varpi_+$ and $\varpi_-$ are minimal representations of the connected group $G_0$.

To fix the notation,
 we suppose $p =1$.  
Then, 
 $G$ is 
the conformal group $O(2,q+1)$ of the 
 Minkowski space
 $\mathbb{R}^{1,q}$, namely, the Euclidean space $\mathbb{R}^{1+q}$
 equipped with the flat Lorentz metric of signature $(1,q)$.
In this case 
our representation $\varpi$ has a long history of study, 
 also in physics
 (see e.g. Todorov \cite{xtodorov}).
The minimal representation $\varpi_+$
may be interpreted as the symmetry of 
the solution space to the mass-zero spin-zero wave equation.
The representation $\varpi_+$ arises also on the Hilbert space of
 bound states of the Hydrogen atom.

\item[5)]
 (spherical case)\enspace
The underlying $({\mathfrak{g}}_{\mathbb{C}}, K)$-module
 $\varpi_K$ has the following $K$-type formula:
\begin{equation}
\label{eqn:K-type}
 \varpi_K \simeq \bigoplus _{\substack{a,b \in \Bbb N, \\ a+\frac p 2 = b + \frac q 2 }}
                     {\mathcal{H}}^a({\mathbb{R}}^{p+1}) \boxtimes  {\mathcal{H}}^b({\mathbb{R}}^{q+1}).  
\end{equation}
In particular,
 the representation $\varpi$ is spherical
 (i.e. contains a non-zero $K$-fixed vector)
 if and only if $p=q$.  

\item[6)]
(infinitesimal character)\enspace
Let ${\mathfrak {Z}}({\mathfrak {g}}_{\mathbb{C}})$ be the center of $U({\mathfrak {g}}_{\mathbb{C}})$.  
Then, 
the infinitesimal character
 of $\varpi_K$ is given by
$$
     (1, \frac{p+q}{2}-1, \frac{p+q}{2}-2, \cdots, 1,0).  
$$
Here,
 we normalize the Harish-Chandra isomorphism for the simple Lie algebra
 ${\mathfrak {g}}_{\mathbb{C}}$ of type $D_n$
 $(n=\frac{p+q}{2}+1)$, 
$$
   \operatorname{Hom}_{{\mathbb{C}}\text{-algebra}}({\mathfrak {Z}}({\mathfrak {g}}_{\mathbb{C}}),{\mathbb{C}})
  \simeq {\mathbb{C}}^n/W(D_n), 
$$
 in a way that the infinitesimal character
 of the trivial one dimensional representation
 is 
$$
     (\frac{p+q}{2}, \frac{p+q}{2}-1, \frac{p+q}{2}-2, \cdots, 1,0).  
$$
\item[7)]
(theta correspondence)\enspace
The representation $\varpi$ is obtained also
 as the theta correspondence of the trivial one-dimensional representation of 
 $SL(2,{\mathbb{R}})$ with respect to the following reductive dual pair
$$
   O(p+1,q+1) \cdot SL(2,{\mathbb{R}}) \subset Sp(p+q+2, {\mathbb{R}}).  
$$
See \cite{xHuZu}.  

\item[8)]
(Gelfand--Kirillov dimension)
\enspace
The Gelfand--Kirillov dimension of the representation $\varpi$,
 to be denoted by $\operatorname{DIM}(\varpi)$, 
 attains its minimum among all unitary representations of $G$, 
 that is, 
$$
   \operatorname{DIM}(\varpi)=p+q-1.  
$$
The associated variety of the underlying $({\mathfrak {g}}_{\mathbb{C}}, K)$-module $\varpi_K$ is given by 
$$
  \operatorname{AV}(\varpi) ={\mathcal{O}}_{\operatorname{min}}^{\mathbb{C}} \cap {\mathfrak {p}}_{\mathbb{C}},
$$
 see \cite[Lemma 4.4]{xkors2}.  
Here, 
 ${\mathcal {O}}_{\operatorname{min}}^{\mathbb{C}}$
 is the minimal nilpotent coadjoint orbit in ${\mathfrak {g}}_{\mathbb{C}}^{\ast}$
 identified with the Lie algebra 
$
     {\mathfrak {g}}_{\mathbb{C}}.  
$

\end{enumerate}

\end{document}